\documentclass[journal]{IEEEtran}

\usepackage{amsmath,amssymb,amsthm,mathrsfs}
\usepackage{mathtools}
\usepackage{graphicx}
\usepackage{algorithm, algpseudocode}
\usepackage{cite}

\newtheorem{theorem}{Theorem}
\newtheorem{lemma}[theorem]{Lemma}
\newtheorem{proposition}[theorem]{Proposition}
\newtheorem{assumption}{Assumption}

\title{Stability Certification for Approximate Nonlinear Model Predictive Control}

\author{An Aloysius Wang, Kostas Margellos, and Mark Cannon%
\thanks{An Aloysius Wang, Kostas Margellos, and Mark Cannon
are with the Department of Engineering Science,
University of Oxford, Oxford OX1 3PJ, UK
(e-mail: an.wang@eng.ox.ac.uk;
kostas.margellos@eng.ox.ac.uk;
mark.cannon@eng.ox.ac.uk). Kostas Margellos is also with the Department of Informatics, National and Kapodistrian University of Athens (e-mail: kmargellos@di.uoa.gr)}}

\begin{document}

\maketitle

\begin{abstract}
Neural-network approximations of nonlinear model predictive control can reduce online computational costs, but ensuring constraint satisfaction and closed-loop stability remains challenging. We propose a framework in which convex optimisation is used to restore feasibility under suitable local assumptions, while bounds on the corrected controller's cost suboptimality provide sufficient conditions for practical closed-loop stability. To assess these conditions over a compact set of initial states, we investigate a branch-and-bound approach based on convex relaxations of the nonlinear dynamics, with probabilistic guarantees established through scenario optimisation. The proposed approach represents a step towards combining the computational efficiency of approximate nonlinear model predictive control with systematic methods for restoring feasibility and assessing stability.
\end{abstract}


\section{Introduction}
Nonlinear model predictive control (NMPC) computes control actions by repeatedly solving a finite-horizon optimisation problem subject to system dynamics and operating constraints. Its ability to accommodate nonlinear dynamics and explicit constraints makes it an attractive approach for a wide range of engineering applications. However, the computational cost of solving a nonlinear optimisation problem at each sampling instant can limit its applicability to systems with stringent real-time requirements. 

Neural-network approximations offer a potential solution by replacing online optimisation with a computationally inexpensive evaluation of a learned control law \cite{Parisini1995}. The practical feasibility of this approach has already been demonstrated experimentally \cite{GomezOrtega1996,Cavagnari1999}, with more recent FPGA-based implementations highlighting its potential for high-speed embedded control \cite{Lucia2021,Chan2021}.

In the context of linear MPC, including time-varying formulations, different neural-network architectures have been investigated with the aim of minimising approximation error \cite{Drgona2018,Kumar2018}. However, small approximation errors alone do not generally guarantee constraint satisfaction or closed-loop stability. To address the former, projection-based architectures have been proposed to restrict neural-network outputs to admissible control inputs defined using control-invariant sets, thereby ensuring recursive feasibility \cite{Chen2018,KargLucia2020}. Separately, primal–dual neural networks have been proposed to assess suboptimality by exploiting weak duality to obtain upper bounds on the optimality gap \cite{Zhang2019}. Neither approach, however, establishes closed-loop stability of the resulting approximate controllers.

One strategy to retain the stability guarantees of MPC is to use neural-network predictions to warm-start an optimisation algorithm \cite{Chen2022}, although the need for an online solve may limit the computational speed-up. Alternatively, in the context of explicit MPC, sufficient conditions for closed-loop stability have been established based on the worst-case approximation error and its Lipschitz constant, both of which can be computed exactly using offline mixed-integer optimisation \cite{FabianiGoulart2023}.

Extensions to nonlinear MPC have proved challenging due to the generally nonconvex nature of the underlying optimisation problem. In this context, feasibility correction has been investigated using a single quadratic programming step with linearised constraints \cite{Vaupel2020}, although satisfaction of the original nonlinear constraints is not guaranteed. Safety-augmentation methods have also been proposed to guarantee constraint satisfaction and closed-loop stability by reverting to a backup trajectory when necessary \cite{Hose2025}, but require a feasible initialisation that may itself involve solving the original NMPC problem.

A complementary approach is to employ robust MPC to account for approximation errors in learned control laws. Statistical bounds on the approximation error have been used to establish probabilistic guarantees for constraint satisfaction and closed-loop stability \cite{Hertneck2018}, while Lipschitz-based bounds have been employed to derive deterministic guarantees under suitable assumptions \cite{Tong2025}. Nevertheless, the associated constraint tightening may introduce conservatism and reduce control performance.

Despite these developments, two important questions remain: how to efficiently restore exact feasibility of neural-network predictions under nonlinear dynamics, and how to certify closed-loop stability of the resulting controller. We address the first using a successive convex optimisation-based correction procedure that asymptotically restores feasibility under suitable local assumptions. For the second, we establish sufficient conditions for practical closed-loop stability based on the cost suboptimality of the corrected controller, and investigate a branch-and-bound method using convex relaxations of the nonlinear dynamics to bound this suboptimality over a compact region of initial states.

The remainder of the paper is organised as follows. Section~\ref{sec:practical stabiltiy} presents the NMPC formulation and establishes sufficient conditions for practical closed-loop stability. Section~\ref{sec:certification} investigates a branch-and-bound method for bounding cost suboptimality using convex relaxations of the nonlinear dynamics, together with scenario optimisation to establish probabilistic guarantees. Section~\ref{sec:numerical} presents numerical results, and Section~\ref{sec:conclusion} concludes the paper. The feasibility-correction procedure is detailed in Appendix~\ref{app:feasibility}, and additional numerical results on computational performance and constraint satisfaction are presented in Appendix~\ref{app:computational_performance}.

\section{Practical Stability of Approximate NMPC}
\label{sec:practical stabiltiy}

Consider a discrete-time nonlinear system
\begin{equation}
    x_{k+1} = f(x_k, u_k),
\end{equation}
where $x_k \in \mathbb{R}^{n_x}$ and $u_k \in \mathbb{R}^{n_u}$ denote the state and control input respectively. The system is subject to state and input constraints
\begin{equation}
    x_k\in\mathcal{X}, \qquad u_k\in\mathcal{U},
\end{equation}
where $\mathcal{X}\subseteq\mathbb{R}^{n_x}$ and $\mathcal{U}\subseteq\mathbb{R}^{n_u}$ are compact sets. We assume that the system admits an equilibrium $(x_{\mathrm{ref}},u_{\mathrm{ref}})$ satisfying $f(x_{\mathrm{ref}},u_{\mathrm{ref}})=x_{\mathrm{ref}},$ which we seek to stabilise.

At each time step $k$, given the current state $x_k$, consider the finite-horizon NMPC problem
\begin{subequations}\label{eq:nmpc}
\begin{align}
    \min_{\mathbf{u}_k}\quad
    &\sum_{j=0}^{N-1}
    \ell(x_{j\mid k},u_{j\mid k})
    +V_f(x_{N\mid k})\\
    \mathrm{s.t.}\quad
    &x_{0\mid k}=x_k,\\
    &x_{j+1\mid k}=f(x_{j\mid k},u_{j\mid k}),
    \quad j=0,\ldots,N-1,\\
    &(x_{j\mid k},u_{j\mid k})
    \in\mathcal{X}\times\mathcal{U},
    \quad j=0,\ldots,N-1,\\
    &x_{N\mid k}\in\mathcal{X}_f.
\end{align}
\end{subequations}
Here, $\mathbf{u}_k=(u_{0\mid k},\ldots,u_{N-1\mid k})$ is the predicted input sequence over the horizon $N$, while $x_{j\mid k}$ denotes the predicted state at stage $j$. Furthermore, $\ell$ is the stage cost, $V_f$ is the terminal cost, and $\mathcal{X}_f\subseteq\mathcal{X}$ is the terminal set. We denote the objective function by $J_N(x_k,\mathbf{u}_k)$ and its optimal value by $V_N(x_k)$. The NMPC problem \eqref{eq:nmpc} is assumed to admit an optimal solution whenever feasible, with the MPC feedback law implemented by applying the first input of an optimal sequence at each time step.

We further impose the following standard assumptions on the stage cost and terminal ingredients.

\begin{assumption}\label{ass:terminal}
There exists a constant $q>0$ such that
\begin{equation}
    \ell(x,u)\geq q\|x-x_{\mathrm{ref}}\|^2
\end{equation}
for all $(x,u)\in\mathcal{X}\times\mathcal{U}$, with $\ell(x_{\mathrm{ref}},u_{\mathrm{ref}})=0$. Furthermore, $V_f$ is nonnegative on $\mathcal{X}_f$, with $x_{\mathrm{ref}}\in\mathcal{X}_f$ and $V_f(x_{\mathrm{ref}})=0$. Moreover, there exists a terminal control law $\kappa_f:\mathcal{X}_f\to\mathcal{U}$ such that, for all $x\in\mathcal{X}_f$,
\begin{align}
    f(x,\kappa_f(x))&\in\mathcal{X}_f,\\
    V_f(f(x,\kappa_f(x)))-V_f(x)
    &\leq-\ell(x,\kappa_f(x)).
\end{align}
\end{assumption}

We now establish sufficient conditions for practical asymptotic stability of a feasible approximate NMPC controller.

\begin{theorem}[Practical asymptotic stability]
\label{thm:practical_stability}
Suppose Assumption~\ref{ass:terminal} holds. Let $c>0$, $\beta\geq 0$, and $0\leq\eta<q$. Define
\begin{equation}
    \Omega_c=\{x:V_N(x)\leq c\}.
\end{equation}
Suppose $\Omega_c$ contains a neighbourhood of $x_{\mathrm{ref}}$ and there exist constants $a,b>0$ such that
\begin{equation}
    a\|x-x_{\mathrm{ref}}\|^2
    \leq V_N(x)
    \leq b\|x-x_{\mathrm{ref}}\|^2
\end{equation}
for all $x\in\Omega_c$.

Furthermore, suppose that for every $x\in\Omega_c$, the approximate controller produces a feasible input sequence $\hat{\mathbf{u}}(x)$ satisfying
\begin{equation}
\label{eq:decrease_condition}
    J_N(x,\hat{\mathbf{u}}(x))-V_N(x)
    \leq \eta\|x-x_{\mathrm{ref}}\|^2+\beta.
\end{equation}

Let
\begin{equation}
    \rho=1-\frac{q-\eta}{b}\in[0,1).
\end{equation}
If
\begin{equation}
    c\geq\beta/(1-\rho),
\end{equation}
then the approximate feedback law $u_k=\hat{u}_{0\mid k}$ is recursively feasible, and $\Omega_c$ is positively invariant. Moreover, for all $x_0\in\Omega_c$ and $k\geq 0$,
\begin{equation}
\label{eq:convergence_rate}
    \|x_k-x_{\mathrm{ref}}\|^2
    \leq 
    \frac{\rho^k}{a}V_N(x_0)+\frac{\beta}{a}
    \frac{1-\rho^k}{1-\rho}.
\end{equation}
In particular,
\begin{equation}
\label{eq:asymptotic_tracking}
    \limsup_{k\to\infty}
    \|x_k-x_{\mathrm{ref}}\|
    \leq\sqrt{\frac{\beta}{a(1-\rho)}}.
\end{equation}
\end{theorem}

\begin{proof}
Suppose $x_k\in\Omega_c$, and let $\hat{\mathbf{u}}_k$ be the feasible input sequence produced by the approximate controller, with corresponding predicted states $\hat{x}_{j\mid k}$. After applying $u_k=\hat{u}_{0\mid k}$, we have $x_{k+1}=\hat{x}_{1\mid k}$.

Consider the shifted input sequence
\begin{equation}
    \tilde{\mathbf{u}}_{k+1}
    =(\hat{u}_{1\mid k},\ldots,
    \hat{u}_{N-1\mid k},
    \kappa_f(\hat{x}_{N\mid k})).
\end{equation}
Under Assumption~\ref{ass:terminal}, $\tilde{\mathbf{u}}_{k+1}$ is feasible at $x_{k+1}$. Moreover,
\begin{align}
    V_N(x_{k+1})
    &\leq J_N(x_{k+1},\tilde{\mathbf{u}}_{k+1})\\
    &\leq J_N(x_k,\hat{\mathbf{u}}_k)
    -\ell(x_k,\hat{u}_{0\mid k}),
\end{align}
where the second inequality follows from the terminal cost decrease condition. Subtracting $V_N(x_k)$ on both sides and applying \eqref{eq:decrease_condition}, we obtain
\begin{equation}
    V_N(x_{k+1})-V_N(x_k)
    \leq -(q-\eta)
    \|x_k-x_{\mathrm{ref}}\|^2+\beta.
\end{equation}

Using $V_N(x_k)\leq b\|x_k-x_{\mathrm{ref}}\|^2$, it follows that
\begin{equation}\label{eq:practical_decrease}
    V_N(x_{k+1})
    \leq \rho V_N(x_k)+\beta,
\end{equation}
where $\rho=1-(q-\eta)/b\in[0,1)$.

Since $c\geq\beta/(1-\rho)$, \eqref{eq:practical_decrease} implies
\begin{equation}
    V_N(x_{k+1})\leq\rho c+\beta\leq c.
\end{equation}
Hence $\Omega_c$ is positively invariant, and recursive feasibility follows by induction.

Iterating \eqref{eq:practical_decrease} gives
\begin{equation}
    V_N(x_k)
    \leq \rho^k V_N(x_0)
    +\beta\frac{1-\rho^k}{1-\rho}.
\end{equation}
Finally, using $a\|x_k-x_{\mathrm{ref}}\|^2\leq V_N(x_k)$, we obtain
\begin{equation}
    \|x_k-x_{\mathrm{ref}}\|^2
    \leq 
    \frac{\rho^k}{a}V_N(x_0)+\frac{\beta}{a}
    \frac{1-\rho^k}{1-\rho}.
\end{equation}
Taking the limit superior as $k\to\infty$ completes the proof.
\end{proof}

Theorem~\ref{thm:practical_stability} establishes that asymptotic stability is recovered when $\beta=0$, while $\beta>0$ yields convergence to a neighbourhood of the reference equilibrium. Crucially, the result requires the approximate input sequence $\hat{\mathbf{u}}(x)$ to satisfy the original nonlinear MPC constraints, which is not generally guaranteed for neural-network predictions. To address this, we present a successive convex optimisation-based correction procedure in Appendix~\ref{app:feasibility} that asymptotically restores feasibility under suitable local assumptions. The remaining challenge is to obtain a uniform bound on the cost suboptimality $J_N(x,\hat{\mathbf{u}}(x))-V_N(x)$ over a region of initial states, particularly since $V_N(x)$ is generally not available in closed form. In the following section, we investigate a branch-and-bound approach based on convex relaxations of the nonlinear dynamics to establish such bounds.

\section{Branch-and-Bound Suboptimality Certification}
\label{sec:certification}
\subsection{Convex Relaxation of the Cost Suboptimality}

Let $\mathcal D\subseteq\mathcal X$ be a nonempty compact polytope of initial states such that $\hat{\mathbf u}(x)$ is feasible for all $x\in\mathcal D$. Throughout this section, we further assume that $\mathcal X$ and $\mathcal U$ are convex polytopes and that a concave overestimator $\overline J_{\mathcal D}$ of $J_N(x,\hat{\mathbf u}(x))$ is available, satisfying
\begin{equation}\label{eq:affine_upper_bound}
    J_N(x,\hat{\mathbf u}(x))
    \leq \overline J_{\mathcal D}(x),
    \qquad \forall x\in\mathcal D.
\end{equation}
For simplicity, we consider the case where $\overline J_{\mathcal D}$ is affine, although more general concave upper bounds could also be used.

For a fixed $0\leq\eta<q$, define
\begin{equation}
r_\eta(x)\coloneqq
J_N(x,\hat{\mathbf{u}}(x))-V_N(x)-\eta\|x-x_{\mathrm{ref}}\|^2,
\end{equation}
and
\begin{equation}\label{eq:beta_star}
\beta_\eta^\star(\mathcal D)
\coloneqq \max\left\{0,\sup_{x\in\mathcal D}r_\eta(x)\right\}.
\end{equation}

Any $\beta\geq\beta_\eta^\star(\mathcal{D})$ satisfies the suboptimality condition of Theorem~\ref{thm:practical_stability} over $\mathcal{D}$.

Let $\mathcal{F}(x)$ denote the set of feasible input sequences for \eqref{eq:nmpc}. Since $V_N(x)=\min_{\mathbf{u}\in\mathcal{F}(x)}
J_N(x,\mathbf{u})$, we may equivalently express
\begin{equation}\label{eq:r_eta_max}
r_\eta(x)=
\max_{\mathbf{u}\in\mathcal{F}(x)}
\bigl[
J_N(x,\hat{\mathbf{u}}(x))-J_N(x,\mathbf{u})
\bigr]-\eta\|x-x_{\mathrm{ref}}\|^2.
\end{equation}
We can therefore obtain an upper bound on $r_\eta(x)$ by replacing the nonlinear dynamics in \eqref{eq:r_eta_max} with convex relaxations, thereby enlarging the feasible set of state–input trajectories.

To construct such a relaxation, let $w=(x,u)$ and suppose that each component of the nonlinear dynamics admits a difference-of-convex (DC) decomposition
\begin{equation}\label{eq:dc_decomposition}
    f_i(w)=g_i(w)-h_i(w),
    \qquad i=1,\ldots,n_x,
\end{equation}
where $g_i$ and $h_i$ are convex on $\mathcal{X}\times\mathcal{U}$. This assumption is relatively mild, since continuous functions on compact convex domains can be approximated arbitrarily well by difference-of-convex functions \cite{HorstThoai1999}.

Let $\{v^{(1)},\ldots,v^{(M)}\}$ denote the vertices of $\mathcal{X}\times\mathcal{U}$. Any $w\in\mathcal{X}\times\mathcal{U}$ can be expressed as
\begin{equation}\label{eq:convex_combination}
    w=\sum_{m=1}^M\lambda_m v^{(m)},
    \quad
    \sum_{m=1}^M\lambda_m=1,\quad \lambda_m\geq 0.
\end{equation}
By convexity of $g_i$ and $h_i$, we have
\begin{equation}\label{eq:dc_bounds}
\begin{aligned}
    g_i(w)-\sum_{m=1}^M\lambda_m h_i(v^{(m)})
    &\leq f_i(w)\\
    &\leq \sum_{m=1}^M\lambda_m g_i(v^{(m)})-h_i(w).
\end{aligned}
\end{equation}

Introducing the predicted states $\mathbf z=(z_0,\ldots,z_N)$ and weights $\boldsymbol\lambda=(\lambda_{j,m})$ as additional
optimisation variables, let $\gamma_\eta(\mathcal D)$ denote the optimal value of
\begin{subequations}\label{eq:relaxed_beta}
\begin{align}
\max_{x,\mathbf z,\mathbf u,\boldsymbol\lambda}\quad
&\overline J_{\mathcal D}(x)
-\eta\|x-x_{\mathrm{ref}}\|^2 \notag\\
&\hspace{2cm}-\sum_{j=0}^{N-1}\ell(z_j,u_j)-V_f(z_N)\\
\mathrm{s.t.}\quad
&x\in\mathcal D,\quad z_0=x,\\
&w_j=\sum_{m=1}^M\lambda_{j,m}v^{(m)},\\
&\sum_{m=1}^M\lambda_{j,m}=1,\quad
\lambda_{j,m}\geq 0,\\
&g_i(w_j)-\sum_{m=1}^M
\lambda_{j,m}h_i(v^{(m)})
\leq z_{j+1,i},\\
&h_i(w_j)+z_{j+1,i}
\leq\sum_{m=1}^M
\lambda_{j,m}g_i(v^{(m)}),\\
&z_N\in\mathcal X_f,
\end{align}
\end{subequations}
where $w_j=(z_j,u_j)$, and the stage constraints are imposed for $j=0,\ldots,N-1$ and $i=1,\ldots,n_x$.

The resulting optimisation problem provides a tractable means of bounding the cost suboptimality of the approximate controller, as formalised in the following proposition.

\begin{proposition}
Under the assumptions of this section,
\begin{equation}\label{eq:beta_upper_bound}
    \beta_\eta^\star(\mathcal D)
    \leq \max\{0,\gamma_\eta(\mathcal D)\}.
\end{equation}
Furthermore, \eqref{eq:relaxed_beta} is a convex optimisation problem provided $\mathcal D$ and $\mathcal X_f$ are convex, and $\ell$ and $V_f$ are convex functions.
\end{proposition}

\begin{proof}
Every feasible trajectory of \eqref{eq:nmpc} satisfies the relaxed dynamics constraints, and $\overline J_{\mathcal D}$ upper bounds the cost of the approximate controller, giving \eqref{eq:beta_upper_bound}.
\end{proof}

The bound in \eqref{eq:beta_upper_bound} may be conservative due to both the relaxation of the nonlinear dynamics and the affine upper bound on $J_N(x,\hat{\mathbf{u}}(x))$. We address these sources of conservatism through two distinct subdivision procedures, beginning with the former while keeping $\mathcal D$ fixed.

\subsection{Refinement of the Dynamics Relaxation}

Let $\mathcal B=(B_0,\ldots,B_{N-1})$, where each $B_j\subseteq\mathcal X\times\mathcal U$ is a polytope restricting the state-input pair at prediction stage $j$. Denote by $\gamma_\eta(\mathcal D;\mathcal B)$ the optimal value of \eqref{eq:relaxed_beta} when the DC relaxation at each stage is constructed using the vertices of $B_j$.

Initially, we take
\begin{equation}
    B_j=
    \begin{cases}
        \mathcal D\times\mathcal U, & j=0,\\
        \mathcal X\times\mathcal U, & j=1,\ldots,N-1.
    \end{cases}
\end{equation}
At each iteration, a node $\mathcal B$ is selected, and one of its domains $B_j$ is partitioned into two polytopes $B_j^{(1)}$ and $B_j^{(2)}$ such that 
\begin{equation}
    B_j=B_j^{(1)}\cup B_j^{(2)}.
\end{equation}
The stage $j$ may be selected based on the largest relaxation gap, evaluated using the optimal weights $\lambda_{j,m}$. This produces two child nodes, with the corresponding DC relaxations constructed over the smaller domains.

Let $\mathscr P$ denote the collection of leaf nodes $\mathcal{B}$ generated by these subdivisions. Since every feasible trajectory of \eqref{eq:nmpc} belongs to at least one node in $\mathscr P$, we have
\begin{equation}\label{eq:bnb_stage_bound}
    \beta_\eta^\star(\mathcal D)
    \leq \hat{\beta}_\eta(\mathcal D)
    \coloneqq \max\left\{0,\,
    \max_{\mathcal B\in\mathscr P}
    \gamma_\eta(\mathcal D;\mathcal B)\right\}.
\end{equation}
For a prescribed tolerance $\beta_{\mathrm{tol}}\geq 0$, any node satisfying $\gamma_\eta(\mathcal D;\mathcal B) \leq\beta_{\mathrm{tol}}$ can be discarded from further subdivision. If all nodes are eventually discarded, then
$\beta_\eta^\star(\mathcal D)\leq\beta_{\mathrm{tol}}$.

The second subdivision procedure addresses the conservatism of the affine upper bound. Suppose that the initial-state region is partitioned into smaller polytopes,
\begin{equation}
    \mathcal D=\bigcup_{s=1}^S\mathcal D_s.
\end{equation}
For each $\mathcal D_s$, a separate affine upper bound $\overline J_{\mathcal D_s}$ can be constructed. Since
\begin{equation}\label{eq:beta_partition}
    \beta_\eta^\star(\mathcal D)
    =\max_{s=1,\ldots,S}
    \beta_\eta^\star(\mathcal D_s),
\end{equation}
the overall bound can be obtained by applying the stage-wise relaxation procedure independently to each $\mathcal D_s$ and taking the maximum of the resulting upper bounds.

It remains to construct the affine upper bounds
$\overline J_{\mathcal D_s}$. We investigate a
sample-based approach to estimating these bounds
and establishing probabilistic guarantees on their
validity.

Let 
\begin{equation}
    G(x)\coloneqq J_N(x,\hat{\mathbf u}(x)),
\end{equation}
and let $\mu_s$ be a probability distribution supported on $\mathcal D_s$. Consider two independent sets of samples, $\{y^{(j)}\}_{j=1}^{M_s}$ and $\{x^{(k)}\}_{k=1}^{K_s}$, drawn independently from $\mu_s$. We solve the convex scenario programme

\begin{subequations}\label{eq:scenario_affine}
\begin{align}
\min_{a_s,b_s}\quad
&\frac{1}{M_s}\sum_{j=1}^{M_s}
\left[
a_s^\top(y^{(j)}-x_c)+b_s-G(y^{(j)})
\right]\\
\mathrm{s.t.}\quad
&G(x^{(k)})\leq
a_s^\top(x^{(k)}-x_c)+b_s,\quad k=1,\ldots,K_s.
\end{align}
\end{subequations}
where $x_c\in\mathcal D_s$. We assume \eqref{eq:scenario_affine} admits a finite optimal solution, denoted by $(a_s^\star,b_s^\star)$, for almost every realisation of the samples. The resulting affine function is
\begin{equation}
    \overline J_{\mathcal D_s}(x)
    \coloneqq (a_s^\star)^\top(x-x_c)+b_s^\star.
\end{equation}

The following proposition establishes a probabilistic guarantee on the cost suboptimality bound obtained using this affine overestimator.

\begin{proposition}
If $\sum_{i=0}^{n_x}\binom{K_s}{i}\varepsilon_s^i(1-\varepsilon_s)^{K_s-i}\leq\delta_s$, then
\begin{equation}
\mathbb P_{\mathcal S_s\sim\mu_s^{M_s+K_s}}
\left\{
\mathbb P_{x\sim\mu_s}
\bigl(r_\eta(x)>\hat{\beta}_\eta(\mathcal D_s)\bigr)
\leq\varepsilon_s
\right\}
\geq 1-\delta_s,
\end{equation}
where $\mathcal S_s$ denotes the collection of fitting and scenario samples.
\end{proposition}

\begin{proof}
Since the fitting samples are independent of the scenario samples, the standard results of convex scenario optimisation \cite{Calafiore2010,Campi2008} apply conditionally on the fitting samples. Averaging over their distribution gives
\begin{equation}\label{eq:scenario_guarantee}
\mathbb P_{\mathcal S_s\sim\mu_s^{M_s+K_s}}
\Bigl\{
\mathbb P_{x\sim\mu_s}
\bigl(
G(x)\leq\overline J_{\mathcal D_s}(x)
\bigr) \geq 1-\varepsilon_s
\Bigr\}
\geq 1-\delta_s.
\end{equation}
Since the DC relaxation provides a valid upper bound whenever $G(x)\leq\overline J_{\mathcal D_s}(x)$, the result follows.
\end{proof}

This provides a probabilistic counterpart to the cost suboptimality condition required by Theorem~\ref{thm:practical_stability}: with confidence at least $1-\delta_s$, the computed upper bound $\hat{\beta}_\eta(\mathcal D_s)$ holds for at least a proportion $1-\varepsilon_s$ of initial states drawn from $\mu_s$. Note, however, that the closed-loop guarantees of Theorem~\ref{thm:practical_stability} require a uniform bound. In particular, if $\Omega_c\subseteq\mathcal D$ and $G(x)\leq\overline J_{\mathcal D}(x)$ for all $x\in\mathcal D$, then, under the remaining conditions of the theorem, robust positive invariance of $\Omega_c$ and the asymptotic tracking-error bound in~\eqref{eq:asymptotic_tracking} follow. Since the probabilistic guarantee allows violations on a subset of states, it does not by itself establish these closed-loop properties.

The overall procedure, combining stage-wise branch-and-bound with subdivision of the initial-state region, is summarised in Algorithm~\ref{alg:certification}.

\begin{algorithm}[t]
\caption{Probabilistic Cost-Suboptimality Certification}
\label{alg:certification}
\begin{algorithmic}[1]
\Require Initial region $\mathcal D$,
tolerance $\beta_{\mathrm{tol}}$,
iteration counts $T,K$,
parameters $\eta,\varepsilon,\delta$
\State $\mathscr Q\gets\{\mathcal D\}$,
$\mathscr C\gets\varnothing$

\For{$t=1,\ldots,T$}
    \If{$\mathscr Q=\varnothing$}
        \State \textbf{break}
    \EndIf

    \State Select $\mathcal D_s\in\mathscr Q$
    and remove it from $\mathscr Q$
    \State Construct $\overline J_{\mathcal D_s}$
    using \eqref{eq:scenario_affine}

    \State Initialise $\mathscr P_s$
    with $\mathcal B=(B_0,\ldots,B_{N-1})$,
    where
    \[
    B_0=\mathcal D_s\times\mathcal U,\quad
    B_j=\mathcal X\times\mathcal U\ (j\geq1)
    \]
    \State Solve the initial convex relaxation

    \For{$k=1,\ldots,K$}
        \State Select $\mathcal B\in\mathscr P_s$
        with the largest upper bound
        \If{$\gamma_\eta(\mathcal D_s;\mathcal B)
        \leq\beta_{\mathrm{tol}}$}
            \State \textbf{break}
        \EndIf

        \State Select stage $j$ using the
        DC relaxation gap
        \State Bisect $B_j$ to obtain
        two child nodes
        \State Solve the child relaxations
        \State Replace $\mathcal B$ in
        $\mathscr P_s$ with its children
    \EndFor

    \State Compute
    $\hat\beta_{\eta,s}\gets
    \max\{0,\max_{\mathcal B\in\mathscr P_s}
    \gamma_\eta(\mathcal D_s;\mathcal B)\}$

    \If{$\hat\beta_{\eta,s}
    \leq\beta_{\mathrm{tol}}$}
        \State Add $\mathcal D_s$ to
        $\mathscr C$
    \Else
        \State Partition $\mathcal D_s$ into smaller subregions
        \State Add the resulting subregions to $\mathscr Q$
    \EndIf
\EndFor

\State \Return $\mathscr C,\mathscr Q$
\end{algorithmic}
\end{algorithm}

\section{Numerical Results}
\label{sec:numerical}

\subsection{Coupled-Tanks system}

We consider a nonlinear coupled-tanks system, with dynamics discretised using the forward Euler method. The state $x_k=(x_{k,1},x_{k,2})^\top \in \mathbb{R}^2$ represents the water levels in the two tanks, while $u_k\in\mathbb R$ denotes the pump input. The discrete-time dynamics
are given by
\begin{equation}\label{eq:tank_dynamics}
f(x,u)=
\begin{bmatrix}
x_1-c_1\sqrt{x_1}+c_u u\\
x_2+c_1\sqrt{x_1}-c_2\sqrt{x_2}
\end{bmatrix},
\end{equation}
where
\begin{equation}
c_1=\frac{\Delta t A_1\sqrt{2g}}{A},
\quad
c_2=\frac{\Delta t A_2\sqrt{2g}}{A},
\quad
c_u=\frac{\Delta t k_p}{A}.
\end{equation}
Here, $A$ denotes the cross-sectional area
of each tank, $A_1$ and $A_2$ are the outlet orifice areas, $k_p$ is the pump coefficient, $g$ is gravitational acceleration, and $\Delta t$ is the sampling period.

The state and input constraints are
\begin{equation}
    \mathcal X=[0.1,30]^2\,\mathrm{cm},
    \qquad
    \mathcal U=[0,24]\,\mathrm{V}.
\end{equation}
We consider regulation about the equilibrium
\begin{equation}
    x_{\mathrm{ref}}=
    \begin{bmatrix}
        17.39645\\15
    \end{bmatrix}\mathrm{cm},
    \qquad
    u_{\mathrm{ref}}=7.27795\,\mathrm{V},
\end{equation}
which satisfies
$f(x_{\mathrm{ref}},u_{\mathrm{ref}})
=x_{\mathrm{ref}}$.

The stage cost is chosen as
\begin{equation}
\ell(x,u)
=(x-x_{\mathrm{ref}})^\top Q(x-x_{\mathrm{ref}})
+R(u-u_{\mathrm{ref}})^2,
\end{equation}
while the terminal cost, terminal controller, and terminal set are given by
\begin{align}
V_f(x)
&=(x-x_{\mathrm{ref}})^\top
P(x-x_{\mathrm{ref}}),\\
\kappa_f(x)
&=u_{\mathrm{ref}}+K_f(x-x_{\mathrm{ref}}),\\
\mathcal X_f
&=\{x:V_f(x)\leq c_f\},
\end{align}
where $P$, $K_f$, and $c_f$ are computed following the procedure in \cite{Doff-Sotta2022}. The system and controller parameters are summarised in Table~\ref{tab:parameters}.

\begin{table}[!t]
\caption{Parameters of the Coupled-Tanks System and NMPC Controller}
\label{tab:parameters}
\centering
\begin{tabular}{@{}lp{0.39\columnwidth}l@{}}
\hline
Symbol & Description & Value (unit) \\
\hline
\multicolumn{3}{c}{Plant parameters} \\
\hline
$A$ & Tank cross-sectional area
    & $15.2\,\mathrm{cm}^2$ \\
$A_1$ & Tank 1 outlet area
    & $0.13\,\mathrm{cm}^2$ \\
$A_2$ & Tank 2 outlet area
    & $0.14\,\mathrm{cm}^2$ \\
$k_p$ & Pump gain
    & $3.3\,\mathrm{cm}^3\mathrm{s}^{-1}\mathrm{V}^{-1}$ \\
$g$ & Gravitational acceleration
    & $981\,\mathrm{cm}/\mathrm{s}^2$ \\
$\Delta t$ & Sampling period
    & $1.4\,\mathrm{s}$ \\
\hline
\multicolumn{3}{c}{NMPC parameters} \\
\hline
$N$ & Prediction horizon & $10$ \\
$Q$ & State cost weight & $I_2\,\mathrm{cm}^{-2}$ \\
$R$ & Input cost weight & $0.2\,\mathrm{V}^{-2}$ \\
$P$ & Terminal cost matrix & $\left[\begin{smallmatrix}
2.1018 & 0.7240\\
0.7240 & 6.6791
\end{smallmatrix}\right]\,\mathrm{cm}^{-2}$ \\
$K_f$ & Terminal feedback gain & $\left[\begin{smallmatrix}
-1.1411 & -0.4189
\end{smallmatrix}\right]\, \mathrm{V}\mathrm{cm}^{-1}$ \\
$c_f$ & Terminal set level & $3.6311$ \\
\hline
\end{tabular}
\end{table}

We approximate the optimal NMPC input sequence using a neural network $\phi_\theta:\mathbb R^2\to\mathbb R^{N}$. The neural-network prediction is defined as
\begin{equation}
\mathbf u_{\mathrm{NN}}(x)
=u_{\mathrm{ref}}\mathbf 1_N
+\phi_\theta(z)
-\exp\left(-\frac{\|z\|^2}{\sigma^2}\right)
\phi_\theta(0),
\end{equation}
where $z = (x-x_\mathrm{ref})/10\, \mathrm{cm}$, $\sigma=0.5,$ and $\mathbf 1_N\in\mathbb R^N$ denotes the vector of ones. This construction ensures that $\mathbf u_{\mathrm{NN}}(x_{\mathrm{ref}}) =u_{\mathrm{ref}}\mathbf 1_N$.

The network $\phi_\theta$ is a fully connected feedforward neural network with one hidden layer of 200 ReLU neurons and a linear output layer of dimension $Nn_u=10$. We generate $10^4$ initial states sampled uniformly from $\mathcal X$, of which 6,372 yield feasible solutions to the NMPC problem using a numerical solver. These state-input pairs form the training dataset. The network is trained by minimising the mean squared prediction error using the Adam optimiser with a learning rate of $10^{-3}$.

Since the neural-network predictions are not necessarily feasible, we apply the successive convex optimisation procedure described in Appendix~\ref{app:feasibility}. On the 6,372 initial states in the dataset, the feasibility rate increases from $38.56\%$ without correction to $99.92\%$ after at most 10 correction iterations. Further results on computational performance and constraint satisfaction are presented in Appendix~\ref{app:computational_performance}.

\subsection{Certification of Practical Stability}

\begin{figure}[!t]
    \centering
    \includegraphics[width=0.8\columnwidth]{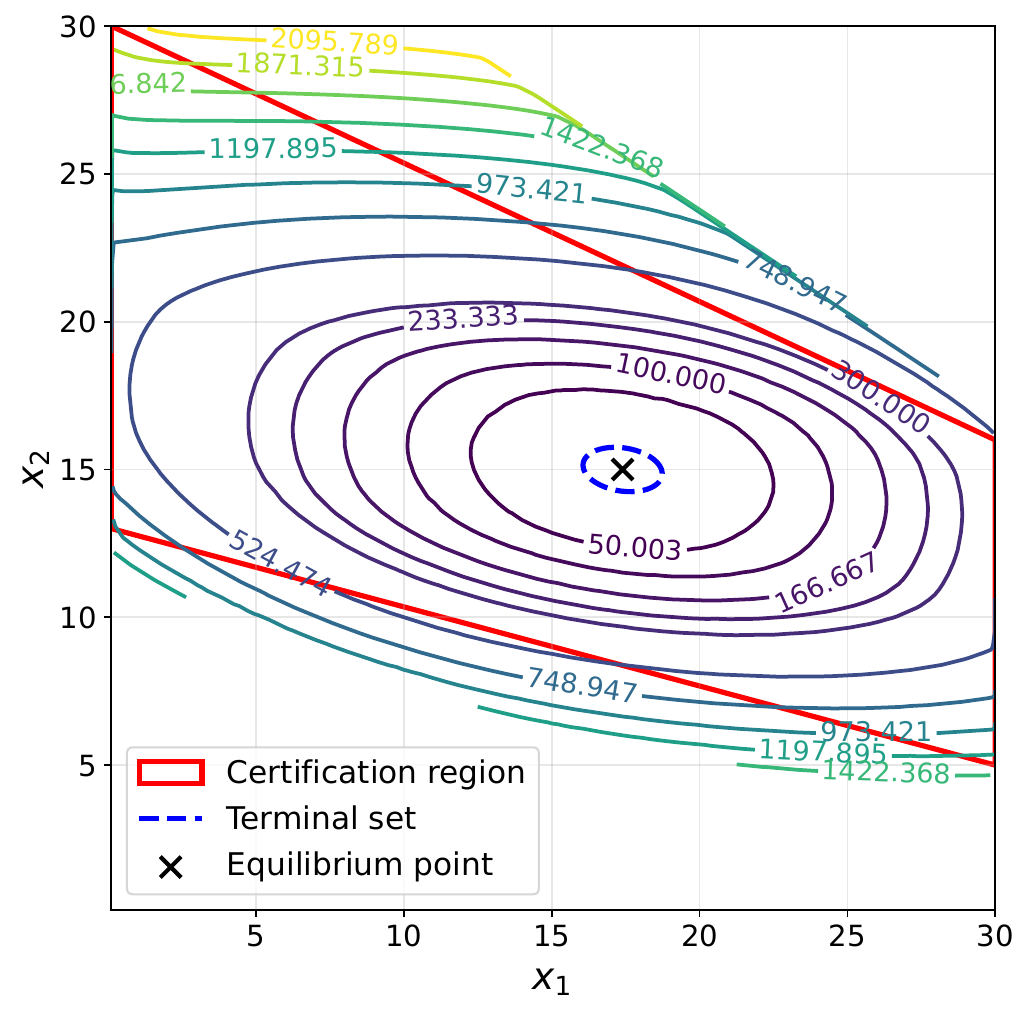}
    \caption{Level sets of the NMPC value function $V_N$, the certification region $\mathcal D$, the terminal set $\mathcal X_f$, and the reference equilibrium $x_{\mathrm{ref}}$.}
    \label{fig:certification_region}
\end{figure}

We then investigate the cost suboptimality of the corrected controller over the initial-state region $\mathcal D\subseteq\mathcal X$, illustrated in Fig.~\ref{fig:certification_region}. The figure also shows the level sets of the NMPC value function $V_N$ and the terminal set $\mathcal X_f$. To this end, we apply Algorithm~\ref{alg:certification} with $T=200$ outer iterations and $K=50$ inner iterations. For each initial-state subregion, an affine upper bound is constructed using 165 scenario samples, with $\varepsilon_s=0.05$ and $\delta_s=0.01$. When subdividing an initial-state subregion, we partition it into four smaller subregions.

\begin{figure}[t]
    \centering
    \includegraphics[width=0.8\columnwidth]{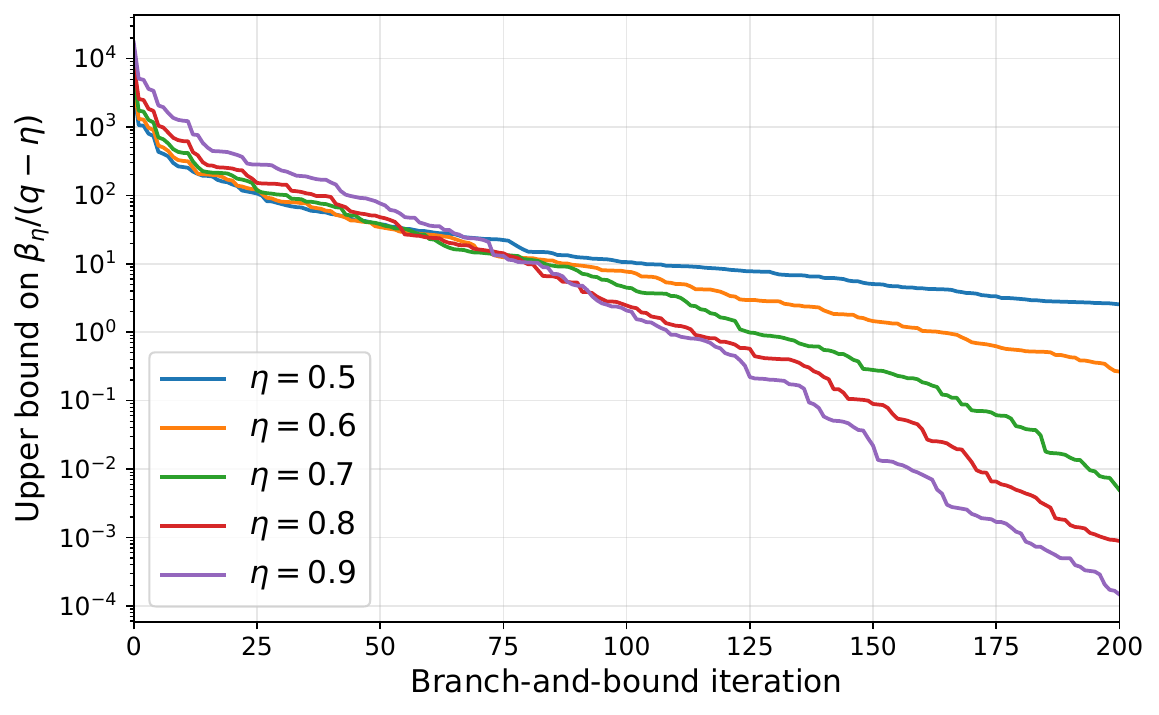}
    \caption{Evolution of the computed upper bound on $\beta_\eta/(q-\eta)$ over successive branch-and-bound iterations for $\eta\in\{0.5,0.6, 0.7,0.8,0.9\}$.}
    \label{fig:beta_convergence}
\end{figure}

Fig.~\ref{fig:beta_convergence} shows the evolution of the computed upper bound on $\beta_\eta/(q-\eta)$ for $\eta\in\{0.5,0.6,0.7,0.8,0.9\}$, with $q=1$. As established in Theorem~\ref{thm:practical_stability}, this quantity is proportional to the upper bound on $\limsup_{k\to\infty}V_N(x_k)$ and therefore determines the practical stability guarantee. The computed bounds decrease with successive branch-and-bound iterations, with larger values of $\eta$ yielding substantially smaller bounds after 200 iterations.

\begin{figure}[t]
    \centering
    \includegraphics[width=0.8\columnwidth]{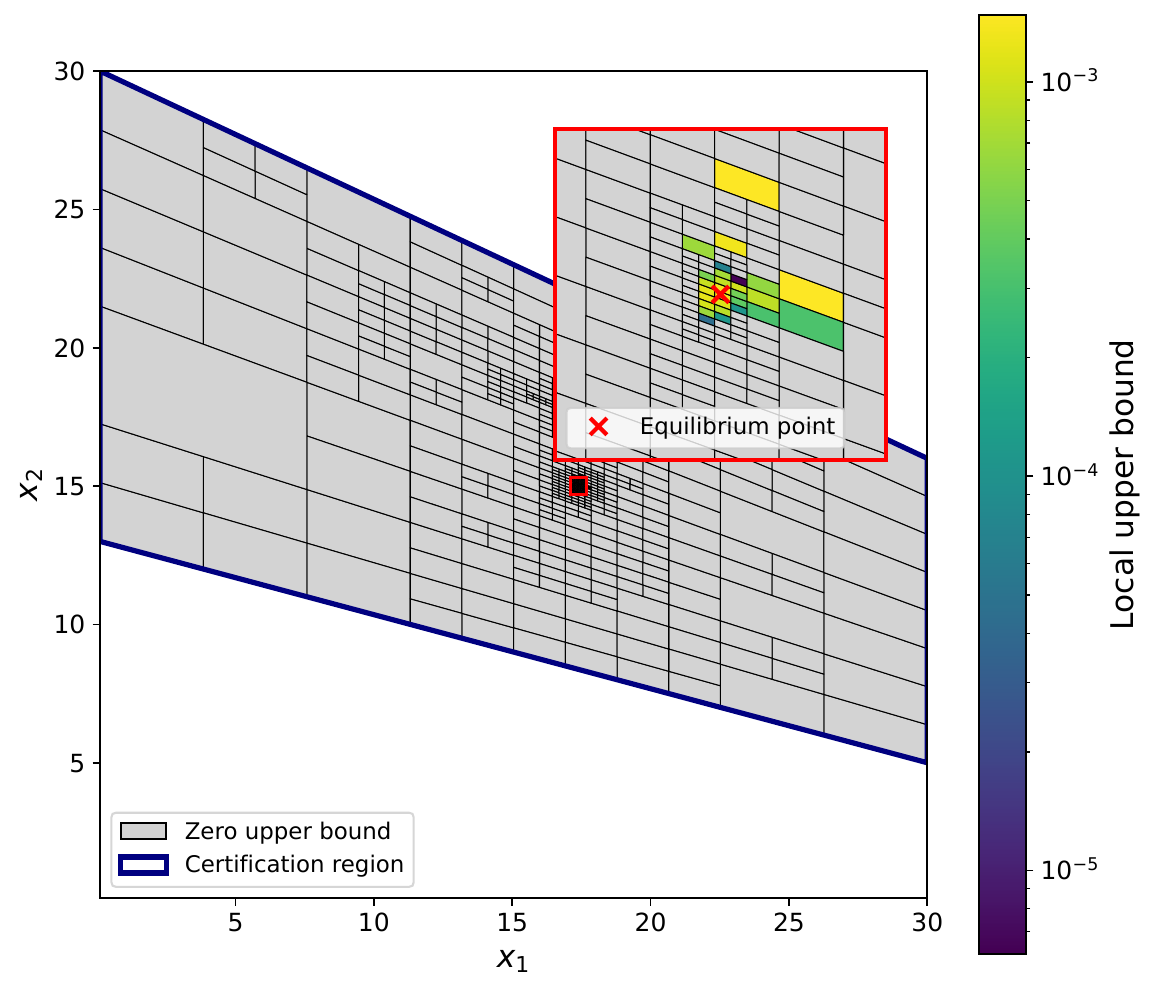}
    \caption{Adaptive subdivision of the certification region after 200 branch-and-bound iterations for $\eta=0.7$. The colour indicates the local upper bound on $\beta_\eta$, with grey corresponding to zero upper bounds. The inset shows the finer subdivisions near the equilibrium $x_{\mathrm{ref}}$.}
    \label{fig:bnb_partition}
\end{figure}

Fig.~\ref{fig:bnb_partition} illustrates the subdivision of the certification region after 200 branch-and-bound iterations for $\eta=0.7$. The colours represent the local upper bounds used to compute $\beta_\eta$. The algorithm produces a nonuniform partition, with finer subdivisions concentrated near the equilibrium, where positive upper bounds remain. In contrast, zero upper bounds are obtained over much of the certification region.

The numerical results suggest that the proposed approach can yield relatively tight practical stability guarantees. In particular, the small values of $\beta_\eta/(q-\eta)$ obtained for larger $\eta$ indicate that the asymptotic bound on the NMPC value function in Theorem~\ref{thm:practical_stability} can be made small. This suggests that the corrected controller can be guaranteed to approach a small neighbourhood of the equilibrium. Note, however, that increasing $\eta$ also increases $\rho$, weakening the convergence rate bound implied by the first term on the right-hand side of~\eqref{eq:convergence_rate}. This highlights a trade-off between the asymptotic bound and the guaranteed convergence rate.

\section{Conclusion}
\label{sec:conclusion}

We have proposed a framework for assessing the stability of neural-network approximations of nonlinear model predictive control. A successive convex optimisation procedure is used to restore feasibility, with local convergence established under suitable regularity assumptions. For the resulting feasible controller, bounds on cost suboptimality provide sufficient conditions for practical closed-loop stability. To assess these bounds over regions of initial states, we developed a branch-and-bound procedure combining convex relaxations of the nonlinear dynamics with affine upper bounds on the corrected controller's cost. The latter are constructed using scenario optimisation, providing probabilistic guarantees on their validity.

\appendices
\section{Feasibility Correction}
\label{app:feasibility}

We consider a successive convex optimisation procedure for correcting an infeasible input sequence. Fix an initial state $x$, and let $\mathbf u\in\mathbb R^{Nn_u}$ denote the stacked input sequence. The corresponding nonlinear state trajectory, obtained by rolling out the dynamics from $x$, is denoted by $X(\mathbf u)\in\mathbb R^{Nn_x}$. 

Let $G(\mathbf u)\coloneqq DX(\mathbf u)$ denote the Jacobian of the nonlinear rollout map $X$. Suppose that the state, input and terminal constraints are polyhedral and can be expressed as
\begin{equation}
    C_xX(\mathbf u)\leq d_x,
    \qquad C_u\mathbf u\leq d_u,
\end{equation}
where $C_x,d_x$ collect the state and terminal constraints over the prediction horizon, while $C_u,d_u$ collect the input constraints. The feasible set of input sequences is then
\begin{equation}
    \mathcal F=
    \{\mathbf u:
    C_xX(\mathbf u)\leq d_x,\;
    C_u\mathbf u\leq d_u\}.
\end{equation}

Given a seed $\mathbf u$, consider the convex optimisation problem
\begin{subequations}\label{eq:correction_qp}
\begin{align}
    \min_{\mathbf x,\mathbf v}\quad
    &\frac12\|\mathbf x-X(\mathbf u)\|^2
    +\frac12\|\mathbf v-\mathbf u\|^2\\
    \mathrm{s.t.}\quad
    &\mathbf x=X(\mathbf u)
    +G(\mathbf u)(\mathbf v-\mathbf u),\\
    &C_x\mathbf x\leq d_x,\qquad
    C_u\mathbf v\leq d_u.
\end{align}
\end{subequations}
Let $T(\mathbf u)=\mathbf v^\ast$ and consider the iteration
\begin{equation}\label{eq:correction_iteration}
    \mathbf u^{j+1}=T(\mathbf u^j).
\end{equation}

We establish local convergence of this procedure to a feasible input sequence. Note that although the analysis below assumes polyhedral constraints, similar results can be established for convex, continuously differentiable constraints under suitable regularity assumptions. Numerical results are presented in Appendix~\ref{app:computational_performance}.

\begin{theorem}\label{thm:correction_convergence}
Suppose that $\mathcal F$ is nonempty, $X$ is twice continuously differentiable, and the linear independence constraint qualification (LICQ) holds at some $\bar{\mathbf u}\in\mathcal F$. Then, for every initial sequence $\mathbf u^0$ sufficiently close to $\bar{\mathbf u}$, \eqref{eq:correction_qp} is feasible at every iteration, and $\mathbf u^j$ converges quadratically to some $\mathbf u^\infty\in\mathcal F$.
\end{theorem}

The proof follows from three lemmas. Define
\begin{equation}
    h(\mathbf u)\coloneqq
    \begin{bmatrix}
        C_xX(\mathbf u)-d_x\\
        C_u\mathbf u-d_u
    \end{bmatrix}.
\end{equation}
Let $\mathcal A$ denote the indices of the active constraints at $\bar{\mathbf u}$. LICQ requires that $Dh_{\mathcal A}(\bar{\mathbf u})$ has full row rank.

\begin{lemma}\label{lem:correction_feasibility}
The convex subproblem \eqref{eq:correction_qp} is feasible in a neighbourhood of $\bar{\mathbf u}$.
\end{lemma}

\begin{proof}
Notice that establishing feasibility of \eqref{eq:correction_qp} requires finding $\mathbf w$ such that
\begin{equation}\label{eq:linearised_h}
    h(\mathbf u)+Dh(\mathbf u)
    (\mathbf w-\mathbf u)\leq0.
\end{equation}

If $\mathcal A$ is empty, the result follows immediately from continuity. Otherwise, by LICQ, $Dh_{\mathcal A}(\bar{\mathbf u})$ has full row rank. Therefore, there exists a direction $p$ such that
\begin{equation}
    Dh_{\mathcal A}(\bar{\mathbf u})p
    =-\mathbf 1.
\end{equation}
By continuity of $Dh$, there exists a neighbourhood $V$ of $\bar{\mathbf u}$ such that
\begin{equation}\label{eq:decreasing_direction}
    Dh_i(\mathbf u)p\leq-\frac12,
    \qquad
    \forall\mathbf u\in V,\quad i\in\mathcal A.
\end{equation}

Define the maximum constraint violation
\begin{equation}
    \rho(\mathbf u)=
    \|[h(\mathbf u)]_+\|_\infty
\end{equation}
and let
\begin{equation}
    \mathbf w(\mathbf u)=
    \mathbf u+2\rho(\mathbf u)p.
\end{equation}
Since $\rho(\mathbf u)\to0$ as $\mathbf u\to\bar{\mathbf u}$, we may choose a smaller neighbourhood $U_0\subseteq V$ such that, for every $\mathbf u\in U_0$, the segment $\{\mathbf u+tp: 0\leq t\leq2\rho(\mathbf u)\}$ lies in $V$. For each $i\in\mathcal A$,
\begin{align}
h_i(\mathbf u)+Dh_i(\mathbf u)
(\mathbf w-\mathbf u) &=h_i(\mathbf u)
+2\rho(\mathbf u)Dh_i(\mathbf u)p\notag\\
&\leq h_i(\mathbf u)-\rho(\mathbf u)
\leq0.
\end{align}
Thus, the active constraints at $\bar{\mathbf u}$ are satisfied.

For any $i\notin\mathcal A$, we have $h_i(\bar{\mathbf u})<0$. By continuity, $h_i(\mathbf u)$ remains strictly negative for $\mathbf u$ sufficiently close to $\bar{\mathbf u}$. Moreover, $\mathbf w-\mathbf u\to0$, and $Dh_i$ is locally bounded. Consequently, sufficiently close to $\bar{\mathbf u}$,
\begin{equation}
    h_i(\mathbf u)+Dh_i(\mathbf u)
    (\mathbf w-\mathbf u)<0.
\end{equation}
By shrinking $U_0$ if necessary, the inactive constraints are satisfied as well. Hence \eqref{eq:linearised_h} holds.
\end{proof}

\begin{lemma}\label{lem:correction_step}
There exist a neighbourhood of $\bar{\mathbf u}$ and a constant $B>0$ such that
\begin{equation}\label{eq:step_bound}
    \|T(\mathbf u)-\mathbf u\|
    \leq B\operatorname{dist}
    (\mathbf u,\mathcal F).
\end{equation}
\end{lemma}

\begin{proof}
Let $U_0$ be as in the previous lemma. By shrinking $U_0$ further, we can ensure that the inactive constraints remain satisfied along the segment $\{\mathbf u+tp: 0\leq t\leq2\rho(\mathbf u)\}$. For $\mathbf u\in U_0$ and $i\in\mathcal A$, we have
\begin{align}
h_i(\mathbf w)
&=h_i(\mathbf u)
+\int_0^{2\rho(\mathbf u)}
Dh_i(\mathbf u+sp)p\,ds\notag\\
&\leq h_i(\mathbf u)-\rho(\mathbf u)
\leq0.
\end{align}
Thus, $\mathbf w\in\mathcal F$. This implies
\begin{equation}\label{eq:distance_violation}
    \operatorname{dist}(\mathbf u,\mathcal F)
    \leq\|\mathbf w-\mathbf u\|
    \leq2\|p\|\rho(\mathbf u).
\end{equation}

Let $\mathbf v=T(\mathbf u)$. By optimality of $\mathbf v$,
\begin{align}
\|\mathbf v-\mathbf u\|^2
&\leq(\mathbf v-\mathbf u)^\top
\bigl(I+G(\mathbf u)^\top G(\mathbf u)\bigr)
(\mathbf v-\mathbf u)\notag\\
&\leq(\mathbf w-\mathbf u)^\top
\bigl(I+G(\mathbf u)^\top G(\mathbf u)\bigr)
(\mathbf w-\mathbf u).
\end{align}
Since $G$ is locally bounded, there exists $M>0$ such that
\begin{equation}
    I+G(\mathbf u)^\top G(\mathbf u)
    \preceq MI
\end{equation}
on $U_0$, after restricting it to a sufficiently small neighbourhood. Hence
\begin{equation}\label{eq:step_rho}
    \|\mathbf v-\mathbf u\|
    \leq2\sqrt M\|p\|\rho(\mathbf u).
\end{equation}

To bound $\rho(\mathbf u)$, choose $B_R(\bar{\mathbf u})\subseteq U_0$ and let $h$ be $L_h$-Lipschitz on $B_{2R}(\bar{\mathbf u})$. Let $\mathbf z$ be a nearest feasible point to $\mathbf u$, which exists since $\mathcal F$ is closed. As $\bar{\mathbf u}\in\mathcal F$,
\begin{equation}
\begin{aligned}
\|\mathbf z-\bar{\mathbf u}\|
&\leq\|\mathbf z-\mathbf u\|
+\|\mathbf u-\bar{\mathbf u}\|\\
&\leq2\|\mathbf u-\bar{\mathbf u}\|
<2R.
\end{aligned}
\end{equation}
Since $h(\mathbf z)\leq0$,
\begin{align}
\rho(\mathbf u)
&\leq\|h(\mathbf u)-h(\mathbf z)\|_\infty
\leq\|h(\mathbf u)-h(\mathbf z)\|
\notag\\
&\leq L_h\|\mathbf u-\mathbf z\|
=L_h\operatorname{dist}
(\mathbf u,\mathcal F).
\end{align}
Combining this with \eqref{eq:step_rho} gives
\begin{equation}
    \|T(\mathbf u)-\mathbf u\|
    \leq B\operatorname{dist}
    (\mathbf u,\mathcal F),
\end{equation}
where $B=2\sqrt M\|p\|L_h$.
\end{proof}

\begin{lemma}\label{lem:correction_quadratic}
There exist a neighbourhood of $\bar{\mathbf u}$ and a constant $K>0$ such that
\begin{equation}\label{eq:quadratic_feasibility}
    \operatorname{dist}
    (T(\mathbf u),\mathcal F)
    \leq K\operatorname{dist}
    (\mathbf u,\mathcal F)^2.
\end{equation}
\end{lemma}

\begin{proof}
Let $R$ be as in the previous lemma, shrunk if necessary so that $B_{(1+B)R}(\bar{\mathbf u}) \subseteq U_0$. For any $\mathbf u\in B_R(\bar{\mathbf u})$, let $\mathbf v=T(\mathbf u)$. Then
\begin{align}
\|\mathbf v-\bar{\mathbf u}\|
&\leq\|\mathbf v-\mathbf u\|
+\|\mathbf u-\bar{\mathbf u}\|\notag\\
&\leq(1+B)\|\mathbf u-\bar{\mathbf u}\|
<(1+B)R.
\end{align}

Let $L_G$ denote the Lipschitz constant of $G$ on $B_{(1+B)R}(\bar{\mathbf u})$, and define $\mathbf x_{\mathrm{lin}} =X(\mathbf u)+G(\mathbf u)(\mathbf v-\mathbf u)$. Then
\begin{equation}
    X(\mathbf v)=\mathbf x_{\mathrm{lin}}+r_T,
    \qquad
    \|r_T\|\leq\frac{L_G}{2}
    \|\mathbf v-\mathbf u\|^2.
\end{equation}
Thus,
\begin{align}
\|[C_xX(\mathbf v)-d_x]_+\|\leq\|C_x\|\|r_T\|
\leq
\frac{L_G\|C_x\|}{2}
\|\mathbf v-\mathbf u\|^2.
\end{align}
Since $C_u\mathbf v\leq d_u$, there is no input constraint violation. Consequently, applying \eqref{eq:distance_violation} at $\mathbf v$, we have
\begin{align}
\operatorname{dist}(\mathbf v,\mathcal F)
\leq2\|p\|\rho(\mathbf v)&\leq L_G\|C_x\|\|p\|
\|\mathbf v-\mathbf u\|^2\notag\\
&\leq K\operatorname{dist}
(\mathbf u,\mathcal F)^2,
\end{align}
where $K=L_G\|C_x\|\|p\|B^2$.
\end{proof}

\begin{proof}[Proof of Theorem~\ref{thm:correction_convergence}]
Suppose $\mathbf u^0$ is sufficiently close to $\bar{\mathbf u}$ such that, with $\delta_j\coloneqq \operatorname{dist}(\mathbf u^j,\mathcal F)$,
\begin{equation}
    K\delta_0\leq\frac12,
    \qquad
    (1+2B)\|\mathbf u^0-\bar{\mathbf u}\|<R.
\end{equation}
We show by induction that $\mathbf u^j\in B_R(\bar{\mathbf u})$ and $\delta_j\leq2^{-j}\delta_0$. Indeed, by Lemma~\ref{lem:correction_quadratic},
\begin{align}
\delta_{j+1}
\leq K\delta_j^2
\leq2^{-j}K\delta_0\delta_j\leq2^{-(j+1)}\delta_j,
\end{align}
and by Lemma~\ref{lem:correction_step}, $\|\mathbf u^{j+1}-\mathbf u^j\|\leq B\delta_j$. Then
\begin{align}
\|\mathbf u^{j+1}-\bar{\mathbf u}\|&\leq\|\mathbf u^0-\bar{\mathbf u}\|
+\sum_{k=0}^j
\|\mathbf u^{k+1}-\mathbf u^k\|\notag\\
&\leq(1+2B)
\|\mathbf u^0-\bar{\mathbf u}\|<R.
\end{align}

Moreover, 
\begin{equation}
    \sum_{j=0}^\infty
    \|\mathbf u^{j+1}-\mathbf u^j\|
    \leq2B\delta_0<\infty
\end{equation}
implies that $\{\mathbf u^j\}$ is Cauchy and converges to some $\mathbf u^\infty$. Since $\delta_j\to0$ and $\mathcal F$ is closed, $\mathbf u^\infty\in\mathcal F$.
\end{proof}

\section{Computational Performance and Constraint Satisfaction}
\label{app:computational_performance}

Table~\ref{tab:computational_performance} compares the mean computation time and maximum constraint violation of the SLSQP-based NMPC solver, the neural-network approximation, and the corrected controller over the test dataset described in Section~\ref{sec:numerical}. The corrected controller achieves an approximately 105-fold reduction in mean computation time compared with SLSQP, while substantially improving constraint satisfaction relative to the uncorrected neural network. These results demonstrate that the feasibility correction can retain much of the computational advantage of neural-network inference while restoring feasibility for the vast majority of tested initial states.

\begin{table}[ht]
    \centering
    \caption{Computational performance and constraint
    satisfaction over the test dataset.}
    \label{tab:computational_performance}
    \begin{tabular}{lcc}
        \hline
        Method & Mean time (ms) & Max. violation \\
        \hline
        NMPC (SLSQP) & 140.2604 & 1.770e-09 \\
        NN & 0.0426 & 2.452e+00 \\
        NN + correction & 1.3321 & 9.903e-06 \\
        \hline
    \end{tabular}
\end{table}

\bibliographystyle{IEEEtran}
\bibliography{references}

\end{document}